\documentstyle[12pt]{article}

\newcommand{\F}{\noindent}

\newcommand{\MP}{\medskip}
\newcommand{\BP}{\bigskip}

\newcommand{\beq}{\begin{eqnarray}}
\newcommand{\ene}{\end{eqnarray}}
\newcommand{\beqn}{\begin{eqnarray*}}
\newcommand{\enen}{\end{eqnarray*}}

\newcommand{\qqqq}{{\mbox{\bf q}}}

\newcommand{\AAA}{{\bf A}}
\newcommand{\BBB}{{\bf B}}

\newcommand{\eq}[1]{(\ref{#1})}

\Large

\begin{document}
\rightline{KIMS-2003-07-07}

\BP

\vskip12pt

\vskip8pt

\begin{center}
\Large
{\bf Does Church-Kleene ordinal $\omega_1^{CK}$ exist?}
\vskip18pt

\normalsize
Hitoshi Kitada
\vskip2pt

Graduate School of Mathematical Sciences

University of Tokyo

Komaba, Meguro-ku, Tokyo 153-8914, Japan

e-mail: kitada@ms.u-tokyo.ac.jp
\vskip10pt
July 7, 2003
\end{center}

\BP

\BP

\leftskip24pt
\rightskip24pt

\small

\noindent
{\it Abstract}: A question is proposed if a nonrecursive ordinal, the so-called Church-Kleene ordinal $\omega_1^{CK}$ really exists.

\leftskip0pt
\rightskip0pt

\vskip 24pt

\large

\large

We consider the systems $S^{(\alpha)}$ defined in \cite{1}.

Let ${\tilde q}(\alpha)$ denote the G\"odel number of Rosser formula or its negation $A_{(\alpha)}$ ($=A_{q^{(\alpha)}}(\qqqq^{(\alpha)})$ or $\neg A_{q^{(\alpha)}}(\qqqq^{(\alpha)})$), if the Rosser formula $A_{q^{(\alpha)}}(\qqqq^{(\alpha)})$ is well-defined.

By ``recursive ordinals" we mean those defined by Rogers \cite{R}. Then that $\alpha$ is a recursive ordinal means that $\alpha<\omega_1^{CK}$, where $\omega_1^{CK}$ is the Church-Kleene ordinal.

\BP

\F
{\bf Lemma}.
The number ${\tilde q}(\alpha)$ is recursively defined for countable recursive ordinals $\alpha<\omega_1^{CK}$. Here `recursively defined' means that ${\tilde q}(\alpha)$ is defined inductively starting from $0$.
\MP

\F
{\bf Remark}. The original meaning of `recursive' is `inductive.' The meaning of the word `recursive' in the following is the one that matches the spirit of Kleene \cite{K} (especially, the spirit of the inductive construction of metamathematical predicates described in section 51 of \cite{K}).
\MP

\F
{\it Proof}. The well-definedness of ${\tilde q}(0)$ is assured by Rosser-G\"odel theorem as explained in \cite{1}.

We make an induction hypothesis that for each $\delta<\alpha$, the G\"odel number ${\tilde q}(\gamma)$ of the formula $A_{(\gamma)}$ ($=A_{q^{(\gamma)}}(\qqqq^{(\gamma)})$ or $\neg A_{q^{(\gamma)}}(\qqqq^{(\gamma)})$) with $\gamma\le\delta$ is recursively defined for $\gamma \le \delta$.

We want to prove that the G\"odel number ${\tilde q}(\gamma)$ is recursively well-defined for $\gamma \le \alpha$.

i) When $\alpha = \delta+1$, by induction hypothesis we can determine recursively whether or not a given formula $A_r$ with G\"odel number $r$ is equal to one of the axiom formulas
$A_{(\gamma)}$ ($\gamma\le\delta$) of $S^{(\alpha)}$. 
In fact, we have only to see, for a finite number
of $\gamma$'s with ${\tilde q}(\gamma)\le r$ and $\gamma\le\delta$, if we have
$A_{(\gamma)}=A_r$ or not. By induction hypothesis that ${\tilde q}(\gamma)$ is recursively well-defined for $\gamma \le \delta$, this is then decided
recursively.

Thus G\"odel predicate $\mbox{\AAA}^{(\alpha)}(a,b)$ and Rosser predicate $\mbox{\BBB}^{(\alpha)}(a,c)$ with
superscript $\alpha$ are recursively defined, and hence are numeralwise
expressible in $S^{(\alpha)}$. Then the Rosser formula $A_{q^{(\alpha)}}(\qqqq^{(\alpha)})$ is well-defined, and the G\"odel number ${\tilde q}(\alpha)$
of Rosser formula or its negation $A_{(\alpha)}$ ($=A_{q^{(\alpha)}}(\qqqq^{(\alpha)})$ or $\neg A_{q^{(\alpha)}}(\qqqq^{(\alpha)})$) is defined recursively.
Thus ${\tilde q}(\gamma)$ is recursively well-defined for $\gamma \le \alpha$.

ii) If $\alpha$ is a countable {\it recursive} limit ordinal, then there is an increasing
sequence of recursive ordinals $\alpha_n<\alpha$ such that 
\beq
\alpha= \bigcup_{n=0}^\infty \alpha_n.\label{al}
\ene
In the system $S^{(\alpha)}$, the totality of the added axioms $A_{(\gamma)}$ $(\gamma<\alpha)$
is the sum of the added axioms $A_{(\gamma)}$ $(\gamma<\alpha_n)$ of $S^{(\alpha_n)}$. By induction hypothesis, ${\tilde q}(\gamma)$ is recursively defined for $\gamma < \alpha_n$. Thus in each $S^{(\alpha_n)}$ we can determine recursively whether or not a given formula $A_r$ is an axiom of $S^{(\alpha_n)}$ by seeing, for a finite number of $\gamma$'s with ${\tilde q}(\gamma) \le r$ and $\gamma<\alpha_n$, if $A_{(\gamma)}=A_r$ or not.

This is extended to $S^{(\alpha)}$. To see this, 
we have only to see the $\gamma$'s with ${\tilde q}(\gamma) \le r$ and
$\gamma<\alpha$, and determine for those finite number of $\gamma$'s if
$A_{(\gamma)}=A_r$ or not. By \eq{al}, 
$$
{\tilde q}(\gamma) \le r\ \mbox{ and }\ \gamma<\alpha
\Leftrightarrow
\exists n \ \mbox{ such that }
{\tilde q}(\gamma) \le r\ \mbox{ and }\ \gamma<\alpha_n.
$$
Then by induction on $n$ with using the result in the above paragraph for $S^{(\alpha_n)}$ and noting that the bound $r$ on ${\tilde q}(\gamma)$ is uniform in $n$, we can show that the condition whether or not ${\tilde q}(\gamma) \le r$ and
$\gamma<\alpha$ is recursively determined. Whence the question whether or not a given formula $A_r$ is one of the axioms $A_{(\gamma)}$ of $S^{(\alpha)}$ with ${\tilde q}(\gamma) \le r$ and $\gamma<\alpha$ is determined recursively.
Thus G\"odel predicate $\mbox{\AAA}^{(\alpha)}(a,b)$ and Rosser predicate $\mbox{\BBB}^{(\alpha)}(a,c)$ with superscript
$\alpha$ are recursively defined, and hence are numeralwise expressible in
$S^{(\alpha)}$. Therefore the Rosser formula $A_{q^{(\alpha)}}(\qqqq^{(\alpha)})$ is well-defined, and the G\"odel number ${\tilde q}(\alpha)$ of Rosser formula
or its negation $A_{(\alpha)}$ ($=A_{q^{(\alpha)}}(\qqqq^{(\alpha)})$ or $\neg A_{q^{(\alpha)}}(\qqqq^{(\alpha)})$) is defined recursively.
Thus ${\tilde q}(\gamma)$ is recursively well-defined for $\gamma \le \alpha$.
This completes the proof of the lemma.

\BP


\F

Assume now that $\alpha$ is a countable limit ordinal such that there is an increasing sequence of recursive ordinals $\alpha_n<\alpha$ with
\beq
\alpha= \bigcup_{n=0}^\infty \alpha_n.\label{al2}
\ene
An actual example of such an $\alpha$ is the Church-Kleene ordinal $\omega_1^{CK}$.

In the system $S^{(\alpha)}$, the totality of the added axioms $A_{(\gamma)}$ $(\gamma<\alpha)$
is the sum of the added axioms $A_{(\gamma)}$ $(\gamma<\alpha_n)$ of $S^{(\alpha_n)}$. By the lemma, ${\tilde q}(\gamma)$ is recursively defined for $\gamma < \alpha_n$. Thus in each $S^{(\alpha_n)}$ we can determine recursively whether or not a given formula $A_r$ is an axiom of $S^{(\alpha_n)}$ by seeing, for a finite number of $\gamma$'s with ${\tilde q}(\gamma) \le r$ and $\gamma<\alpha_n$, if $A_{(\gamma)}=A_r$ or not.

This is extended to $S^{(\alpha)}$. To see this, 
we have only to see the $\gamma$'s with ${\tilde q}(\gamma) \le r$ and
$\gamma<\alpha$, and determine for those finite number of $\gamma$'s if
$A_{(\gamma)}=A_r$ or not. By \eq{al2}, 
$$
{\tilde q}(\gamma) \le r\ \mbox{ and }\ \gamma<\alpha
\Leftrightarrow
\exists n \ \mbox{ such that }
{\tilde q}(\gamma) \le r\ \mbox{ and }\ \gamma<\alpha_n.
$$
Then by induction on $n$ with using the above result for $S^{(\alpha_n)}$ in the preceding paragraph and noting that the bound $r$ on ${\tilde q}(\gamma)$ is uniform in $n$, we can show that the condition whether or not ${\tilde q}(\gamma) \le r$ and
$\gamma<\alpha$ is recursively determined. Then within those finite number of $\gamma$'s with ${\tilde q}(\gamma)\le r$ and $\gamma<\alpha$, we can decide recursively if for some $\gamma<\alpha$ with ${\tilde q}(\gamma)\le r$, we have $A_r=A_{(\gamma)}$ or not. Therefore we can determine recursively whether or not a given formula $A_r$ is an axiom of $S^{(\alpha)}$.

Therefore G\"odel predicate $\mbox{\AAA}^{(\alpha)}(a,b)$ and Rosser predicate $\mbox{\BBB}^{(\alpha)}(a,c)$ are recursively defined, and hence are numeralwise expressible in
$S^{(\alpha)}$. Then the G\"odel number $q^{(\alpha)}$ of the formula
$$
\forall b [\neg A^{(\alpha)}(a,b)\vee \exists c(c\le b \hskip3pt\&\hskip2pt B^{(\alpha)}(a,c))]
$$
is well-defined, and hence Rosser formula $A_{q^{(\alpha)}}(\qqqq^{(\alpha)})$ is well-defined and Rosser-G\"odel theorem applies to the system $S^{(\alpha)}$. Therefore we can extend $S^{(\alpha)}$ consistently by adding one of Rosser formula or its negation $A_{(\alpha)}$ ($=A_{q^{(\alpha)}}(\qqqq^{(\alpha)})$ or $\neg A_{q^{(\alpha)}}(\qqqq^{(\alpha)})$) to the axioms of $S^{(\alpha)}$ and get a consistent system $S^{(\alpha+1)}$.

\BP

In particular if we assume a least nonrecursive ordinal $\omega_1^{CK}$ exists and take $\alpha=\omega_1^{CK}$, we get a consistent system $S^{(\omega_1^{CK}+1)}$. This contradicts the case ii) of the theorem in \cite{1}. We now arrive at

\BP

\F
{\bf Question}. The least nonrecursive ordinal, the so-called Church-Kleene ordinal $\omega_1^{CK}$ has been assumed to give a bound on recursive construction of formal systems (see \cite{F}, \cite{S}, \cite{T}). However the above argument seems to question if $\omega_1^{CK}$ really exists in usual set theoretic sense. How should we think?

\BP

\end{document}